
\input amstex
\documentstyle{amsppt}
\magnification=1200
\catcode`\@=11
\redefine\logo@{}
\catcode`\@=13
\define\wed{\wedge}
\define \bn{\Bbb N}
\define \bz{\Bbb Z}
\define \bq{\Bbb Q}
\define \br{\Bbb R}
\define \bc{\Bbb C}
\define \bh{\Bbb H}
\define \bp{\Bbb P}
\define \ot{\bp\bh_t^{+}}
\define \gt{\Gamma_t}
\define \gts{\Gamma^*_t}
\define \gti{\widetilde\Gamma_t}
\define \km{\hbox{Km}}

\TagsOnRight
\NoBlackBoxes

\topmatter

\title
Minimal Siegel modular threefolds
\endtitle

\author
Valeri Gritsenko
\footnote{Supported by  Institute Fourier, Grenoble,
the  Schwerpunktprogramm
``Komplexe Mannigfaltigkeiten''  grant Hu 337/2-4 and DFG grant 436 Rus
17/108/95.
\hfill\hfill}
and
Klaus Hulek
\endauthor

\address
St. Petersburg Department Steklov Mathematical Institute
Fontanka 27,
\newline
${}\hskip 9pt $
191011 St. Petersburg,  Russia
\endaddress

\email
gritsenk\@cfgauss.uni-math.gwdg.de;\ \
gritsenk\@pdmi.ras.ru
\endemail

\address
Institut f\"ur Mathematik, Universit\"at Hannover,
Postfach 6009,
\newline
${}\hskip 9pt $
D-30060 Hannover, Germany
\endaddress

\email
hulek\@math.uni-hannover.de
\endemail

\abstract
In this paper we study the maximal extension $\Gamma_t^*$ of the subgroup
$\Gamma_t$ of $\operatorname{Sp}_4
(\bq)$ which is conjugate to the paramodular group.
The index of this extension is $2^{\nu(t)}$ where $\nu(t)$ is the
number of prime divisors of $t$. The group $\Gamma_t^*$ defines the
minimal modular threefold ${\Cal A}_t^*$ which is a finite quotient of
the moduli space ${\Cal A}_t$ of $(1,t)$-polarized abelian surfaces.
A certain degree 2 quotient of
${\Cal A}_t$ is a moduli space of lattice polarized $K3$ surfaces.
The space ${\Cal A}_t^*$ can be interpreted as the space of Kummer
surfaces associated to $(1,t)$-polarized abelian surfaces.
Using the action of $\Gamma_t^*$
on the space of Jacobi forms we show that many spaces between ${\Cal A}_t$
and ${\Cal A}_t^*$ posess a non-trivial 3-form, i.e. the Kodaira
dimension of these spaces is non-negative. Finally we determine the
divisorial part of the ramification locus of the finite map ${\Cal
A}_t\rightarrow {\Cal A}_t^*$ which is a union of Humbert surfaces.
We interprete the corresponding
Humbert surfaces as Hilbert modular surfaces.
\endabstract
\rightheadtext
{Minimal Siegel modular threefolds}

\leftheadtext{V. Gritsenko and  K. Hulek}
\endtopmatter

\document

\document

\head
Introduction
\endhead

The moduli space ${\Cal A}_t$ of abelian surfaces with a
$(1,t)$-polarization is the quotient of the Siegel upper half plane
$\bh_2$ by a subgroup $\Gamma_t$ of
\,$\operatorname{Sp}_4(\bq)$ which is conjugate to the
paramodular group ${\widetilde \Gamma}_t$.
In \S\ 1 we define an
isomorphism between the symplectic group  and the special
orthogonal group $\operatorname{SO}(3,2)$ over the integers.

This exhibits
$\Gamma_t/\{\pm E_4\}$ as a subgroup of the orthogonal group
$\operatorname{SO}(L_t)$ where
$L_t$ is the lattice of rank 5 equipped with the form
$<2t>\hskip-1pt\oplus\, 2\,U$
(here $U$ denotes the hyperbolic plane).  Let ${\widehat L}_t$ be the
dual lattice of
$L_t$. The image of $\Gamma_t$ in $\operatorname{O}(L_t)$ acts trivially
on ${\widehat L}_t/L_t$. The orthogonal group
$\operatorname{O}({\widehat L}_t/L_t)$ is isomorphic to
$(\bz/2\bz)^{\nu(t)}$ where ${\nu}(t)$ is the number of prime divisors of
$t$. For every $d||t$ (i.e. $d|t$ and $(d,t/d)=1)$ we construct an
element $V_d$ in $\operatorname{Sp}_4(\br)$.
These elements $V_d$ define a normal
extension $\Gamma_t^*$ of $\Gamma_t$ of index $2^{\nu(t)}$ such that
$\Gamma_t^*/ \Gamma_t\cong \operatorname{O}({\widehat L}_t/L_t)$.
It turns out that
$\Gamma_t^*$ is the maximal normal extension of $\Gamma_t$ as
a discrete  subgroup of $\operatorname{Sp}_4(\br)$.
Hence we can consider the moduli space
$\Cal A_t^*=\Gamma_t^*\backslash\bh_2$ as ``minimal'' Siegel modular
threefolds. In $\S\ 1$ we also give a geometric interpretation of the
action of $V_d$ on the moduli space
${\Cal A}_t$. In particular $V_t$ identifies a polarized abelian surface
with its dual. It also turns out that the space
$(\Gamma_t \cup \Gamma_tV_t)\backslash \bh_2$
is isomorphic to the moduli space of
lattice polarized $K3$ surfaces with a polarization of type
$<2t>\hskip-1pt\oplus\, 2E_8(-1)$.
Lattice polarized $K3$ surfaces have been studied by
Nikulin \cite{N2}.
They play a role in mirror symmetry for $K3$ surfaces
(see  Dolgachev \cite{D}).
Moreover the variety ${\Cal A}_t^*$ is the space of Kummer surfaces
associated to $(1,t)$-polarized abelian surfaces.
(For a precise statement see Theorem 1.5.)

In $\S\ 2$ we study the action of the elements $V_d$ on the space
of Jacobi forms. This gives rise to a decomposition of the space of Jacobi
forms which was originally found by Eichler and Zagier \cite{EZ}.
Using lifting results due to the first author this enables us
to prove that many moduli spaces lying between ${\Cal A}_t$ and
the minimal Siegel modular threefold ${\Cal A}_t^*$ are not unirational,
resp. have non-negative Kodaira dimension.
This method, however, unfortunately does not give us
information about the Kodaira dimension of ${\Cal A}_t^*$ itself.

If one wants to determine the Kodaira dimension of ${\Cal A}_t^*$ one
needs precise information about the ramification locus of the finite map
${\Cal A}_t\rightarrow {\Cal A}_t^*$.
This turns out to be a difficult problem.
In $\S\ 3$ we determine the divisorial part of this ramification locus for
square free $t$ (the general case can be treated by the same method). The
divisorial part of this ramification locus is a finite union of Humbert
surfaces.
To determine these surfaces we reexamine the theory of Humbert
surfaces from the point of view of the orthogonal group.
This turns out to be a very useful way of studying Humbert surfaces.
An example, originally
due to Brasch, shows that the ramification locus can also contain curve
components. We finally interprete the Humbert surfaces in the ramification
locus as Hilbert modular surfaces.

\vskip 0.5truecm
\head
\S\ 1. The symplectic and  orthogonal groups
\endhead

The local isomorphism between the symplectic group
$\operatorname{Sp }_4(\Bbb R)$
and  the special orthogonal
group $\operatorname{SO }(3,2)_\Bbb R$ of signature $(3,2)$
is well known.
In this section we define this  isomorphism  over $\Bbb Z$.

Let us fix a lattice
$$
L=e_1\Bbb Z\oplus e_2\Bbb Z\oplus e_3\Bbb Z\oplus e_4 \Bbb Z.
$$
We identify $l\in L$ with a column-vector in the basis $\{e_i\}$.
$L^2=L\wedge L$ is the lattice of integral bivectors,
which is isomorphic to the lattice of integral skew-symmetric matrices.
The bivector $e_i\wedge e_j$ corresponds to the  elementary
skew-symmetric matrix
$E_{ij}$, which has only two non-zero elements $e_{ij}=1$ and
$e_{ji}=-1$.
Any linear transformation $g:\,L\to L$ induces a linear map
$
\wedge^2 g: L\wedge L\to L\wedge L
$
on the $\Bbb Z$-lattice of bivectors.
If $g$ is  represented  with respect to the  basis $\{e_i\}$ by
the matrix $G$, then
$$
(\wedge^2 g)\,(X)=G X \,{}^t G
\qquad\text{for any }\  X=\sum_{i<j}x_{ij}e_i\wedge e_j\in L\wedge L.
$$
One can define a symmetric bilinear form $(X,Y)$ on   $L\wedge L$
$$
X\wedge Y=(X,Y)\,e_1\wedge e_2\wedge e_3\wedge e_4\in \wedge^4 L.
$$
It is known, that
$(X, X)=2\, \hbox{Pf\,} (X)$,
where $\hbox{Pf\,}(X)$ is the  Pfaffian of the matrix  $X$, and
$\hbox{Pf\,}(MX{}^tM)=\hbox{Pf\,}(X)\, \hbox{det\,} M$.
\definition{Definition}The group
$$
\gti=\{g:L\to L\,|\ \wedge^2 g\,(W_t)=W_t,\quad\text{where }
\ W_t=e_1\wedge e_3+t e_2\wedge e_4\}                    \tag{1.1}
$$
is called the integral {\it paramodular group} of level $t$.
\enddefinition

The lattice $L_t=W_t^{\perp}$ consisting  of all elements of $L\wedge L$
orthogonal to $W_t$ has the following  basis
$$
L_t=(e_1\wedge e_2,\, e_2\wedge e_3,\, e_1\wedge e_3 -te_2\wedge e_4,\,
e_4\wedge e_1,\, e_4\wedge e_3)\,\Bbb Z^5.
$$
We fix this basis for the rest of the paper.
The symmetric bilinear form $(\cdot, \cdot)$  defines a  quadratic form
$S$ of signature $(3,2)$ on the lattice $L_t$, which has
the following form in the given basis
$$
S_t=\pmatrix
0&0&0&0&-1\\
0&0&0&-1&0\\
0&0&2t&0&0\\
0&-1&0&0&0\\
-1&0&0&0&0
\endpmatrix.   \tag{1.2}
$$

The group of the real points of the paramodular group is conjugate
to $\operatorname{Sp }_4(\br)$. Thus the determinant  of any element of
the paramodular group  equals  one and   $\wed^2 g$ keeps
the bilinear form on $L\wedge L$.

This  gives us a homomorphism from the symplectic group in the orthogonal
group of the isometries of  the lattice $L_t$
$$
\wed^2:\, \gti \to \operatorname{O}(L_t).
$$
The paramodular group $\gti$ is conjugate to a
subgroup of the usual rational symplectic group:
$$
\gt:=I_t^{-1}\gti I_t=\left\{\pmatrix
*    &   *   &   *   &   t*\\
t*   &   *   &  t*   &   t*\\
{*}   &   *   &   *   &   t*\\
{*}  &  t^{-1}*& *   &   *
\endpmatrix
 \in \operatorname{Sp }_4 (\bq) \right\},
$$
where all entries $*$ denote integers and
$I_t=\hbox{diag }(1,1,1,t)$.

The quotient space
$$
{\Cal A}_t=\gt\backslash \bh_2
$$
is the coarse moduli space of abelian surfaces with a polarization of
type $(1,t)$.

The composition of the conjugation
with the homomorphism $\wed^2$ defines a homomorphism
$$
\Psi: \gt\to \operatorname{O}(L_t)\qquad\text{where }\
\Psi(g)=\wed^2(I_t\, g \,I_t^{-1}).
\tag{1.3}
$$
One can extend $\Psi$ to  the real symplectic group
$\Gamma_t(\br)\cong \operatorname{Sp }_4(\br)$.

Let
$
\widehat L_t=\{u\in L_t\otimes \Bbb Q\,|\
\forall\ l\in  L_t\  (l,u)\in \Bbb Z\}
$
be the dual lattice of $L_t$. The discriminant group
$$
A_{t}:=\widehat L_t/L_t = (2t)^{-1} \bz/\bz \cong \bz/2t \bz
$$
is a finite abelian group equipped with a quadratic form
$$
q_t:\, A_t\times A_t\to (2t)^{-1} \bz/2\bz\qquad\qquad
q_t(l,l)\equiv (l,l)_{\widehat L_t}\text{\,mod\,}2\bz
$$
(see \cite{N1} for a general definition).
Any $g\in \operatorname{O}(L_t)$ acts
on the  finite group $A_t$.
By
$$
\widehat{\hbox{O}}(L_t)=\{g\in \operatorname{O}(L_t)\,|\
\forall \ell\in \widehat L_t\quad g\ell-\ell\in L_t\}
$$
we denote the subgroup of the orthogonal group consisting of elements
which act identically on the discriminant group.

One can easily prove the next lemma (see \cite{G1}).

\proclaim{Lemma 1.1}The following relations are valid

{\rm 1}. $\Psi(\gt)\subset \widehat{\hbox{\rm SO}}(L_t)
=\widehat{\hbox{\rm O}}(L_t)\cap \operatorname{SO }(L_t)$;

\smallskip
{\rm 2}. $\hbox{\rm Ker}\,\Psi=\{\pm E_4\}$.
\endproclaim

The finite orthogonal group $\operatorname{O}(A_t)$
can be described as follows.
For every $d||t$ (i.e. $d|t$ and $(d,\dsize\frac{t}d)=1$)
there exists a unique  ($\operatorname{mod\,} 2t$) integer $\xi_d$
satisfying
$$
\xi_d=-1\, \operatorname{mod }\, 2 d,\quad \xi_d=
1\, \operatorname{mod }\, 2t/d.
$$
All such  $\xi_d$ form the group
$$
\Xi (t)=\{\,\xi \operatorname{mod} 2 t\ |\
\xi^2= 1 \operatorname{mod }4 t\,\}
\cong (\bz/2\bz)^{\nu (t)},                       \tag{1.4}
$$
where $\nu(t)$ is the number of prime divisors of $t$. It is evident
that
$\operatorname{O}(A_t)\cong \Xi (t)$.

One can take an element in $\operatorname{SO }(L_t)$ realising the
multiplication by
$\xi_d$ on $A_t$. It  gives us an element  in $\operatorname{Sp }_4(\br)$
with  integral $\Psi$-image. For example,
for every  $d||t$ we can define $x,y \in \bz$ (which are
not uniquely determined) such that
$$
xd-yt_d=1\quad \operatorname{where }\  t_d=\frac t d.
$$
The matrix
$$
\widetilde{V}_d=\pmatrix
dx   &   -1   &   0   &   0\\
-yt   &    d   &   0   &   0\\
0    &    0   &   d   &   yt\\
0    &    0   &   1   &   dx
\endpmatrix
$$
is an integral symplectic similitude of degree $d$.
We put
$$
V_d=\frac 1{\sqrt d} {\widetilde V}_d \in \operatorname{Sp }_4 (\br).
$$
$V_d$ has the following $\Psi$-image
$$
\Psi({V}_d)=
\pmatrix
1 &  0&   0   &   0  &0\\
0 &d  & -2yt  &y^2t_d&0\\
0 &-1 &dx+t_dy& -xy  &0\\
0 &t_d& -2tx  &x^2d  &0\\
0 & 0 &   0   &  0   &1
\endpmatrix.                                              \tag{1.5}
$$
We note here that
$$
xd+t_dy= -1\, \operatorname{mod }\, 2d
\qquad\text{and}\qquad xd+t_dy= 1
\, \operatorname{mod }\, 2t_d,
$$
thus $V_d$ induces the multiplication by $\xi_d$ on $A_t$.

It is easy to see, that for all  $V_d$ ($d||t$)
$$
V_d^2\in \gt,  \qquad          V_d\gt V_d=\gt.
$$
I.e. $V_d$ are involutions modulo $\gt$.
Therefore
one can define the following normal extension of the paramodular group
$\gt$.
\definition{Definition}
$\gts$   is  the  group generated by the  elements of  $\gt$ and  $V_d$
for all $d||t$.
\enddefinition

In accordance with  Lemma 1.1 any element
in  $\Psi(V_d\gt)$ defines
the same automorphism of $A_t$,
thus
$$
\gts/\gt\cong \operatorname{O}(\widehat L_t/L_t)\cong \Xi (t)
\cong (\bz/ 2\bz)^{\nu(t)}.         \tag{1.6}
$$
\medskip

The real orthogonal group
$\operatorname{O}_{\br}(L_t)=\operatorname{O}(L_t\otimes \br)$
acts on  a domain
lying  on a projective quadric, more exactly on
$$
\bp\bh_t^{\,3}=\bp\bh_{\,L_t}^{\,3}=
\{Z\in \bp(L_t\otimes \bc)\,|\ (Z,Z)=0,\ (Z,\overline Z)<0\}=
\bp{\bh}_t^{+} \cup \overline{\bp\bh}_t^{+},
$$
where
$$
{\bp\bh}_t^{+}=
\{Z={}^t\bigl((t z_2^2-z_1z_3),\,z_3,\,z_2,\,z_1,\,1\bigr)
\cdot z_0\in \bp\bh_t^{3}\
|\ \hbox{Im}\,(z_1)>0\}.   \tag{1.7}
$$
This  is a classical homogeneuos domain of type IV.
The condition $(Z,\overline Z)<0$ is equivalent to
$$
y_1y_3-ty_2^2>0
\qquad\text{where }\ y_i=\hbox{Im}\,(z_i).
$$
Taking $z_0=1$ one  gets the corresponding cylindric domain in the affine
coordinates $(z_i)_{1\le i\le 3}$
$$
{\bh}_t^{+}=
\{Z={}^t\bigl(z_3,\,z_2,\,z_1)\in \bc^3\,|\
\ y_1y_3-ty_2^2>0,\  \ y_1=\hbox{Im}\,(z_1)>0\}.
$$
The domain ${\bh}_t^{+}$ for $t=1$ coincides with
the Siegel upper half-plane $\bh_2$.
For a  general $t$ one can define the  following isomorphism
of the  complex domains
$$
\psi_t: \bh_2\to \bh_t^+\qquad
\psi_t(\pmatrix \tau_1&\tau_2\\\tau_2&\tau_3\endpmatrix)=
{}^t\bigl(\frac{\tau_3}t,\, \frac{\tau_2}t,\, \tau_1\bigr). \tag{1.8}
$$

The linear action of the real orthogonal group
$\operatorname{O}_\br(L_t)$
on $\bp\bh_t^{+}$ defines  ``fractional-linear'' transformations
on $\bh_t^+$.
By $\operatorname{O}_\br^+(L_t)$ we denote the  subgroup  of
index $2$ of
the orthogonal group
consisting  of elements  which leave $\bp\bh_t^+ $ invariant.
(This is the  subgroup of the  elements with
real spin norm equal one.)

\proclaim{Proposition 1.2}Let $t$ be square free.
Then  $\Psi$ defines the following  isomorphisms
$$
\Psi:\ \gts/\{\pm E_4\} \to \operatorname{SO }^+(L_t),
$$
where $\operatorname{SO }^+(L_t)
=\operatorname{SO}(L_t)\cap \operatorname{O }^+_\br(L_t)$, and
$$
\Psi:\ \gt/\{\pm E_4\} \to \widehat{\hbox{\rm SO}}^+(L_t),
$$
where $\widehat{\hbox{\rm SO}}^+(L_t)=\operatorname{SO }^+
(L_t)\cap \widehat{\hbox{\rm SO}}(L_t)$.
Moreover  the following diagram is commutative
$$
\CD
\bh_2@>g>>\bh_2\\
@V\psi_t VV @V\psi_t  VV\\
\bh_t^+@>\Psi(g)>>\bh_t^+.
\endCD
$$
\endproclaim
\demo{Proof}The diagram is commutative for any
$g\in \operatorname{Sp }_4(\br)$.
To prove this one has to  calculate
the images   of  standard generators of $\operatorname{Sp }_4(\br)$
under $\Psi$.

It is known  that for square free $t$,
the group $\hbox{P}\gts\cong \gts/\{\pm E_4\}$ is
a maximal discrete subgroup of the group of  analytic
automorphisms of  $\bh_2$ and
$[\gts : \gt]=2^{\nu(t)}$
(see for example \cite{Al}, \cite{Gu}).
{}From the description of the finite orthogonal group
$\operatorname{O}(A_t)$ given
in \thetag{1.4} we obtain that
$[\operatorname{SO }^+(L_t):\widehat{\hbox{SO}}^+\hskip -2pt(L_t)]
=2^{\nu(t)}$.
The statement of the proposition about the  isomorphism of the groups
follows
from the maximality of $\gts$ and Lemma 1.1.
\newline\qed\enddemo

The coset $V_t\,\gt$  (in the  case $d=t$ we may take $x=0, y=-1$)
can also be written in the form
$$
V_t\,\gt=
\pmatrix
0   &   \sqrt{t}^{-1}   &   0   &   0\\
\sqrt{t}  & 0 &   0   &   0\\
0   &   0   &   0   &   \sqrt{t}\\
0   &   0   &   \sqrt{t}^{-1}   &   0
\endpmatrix
\gt.
$$
According to \thetag{1.3} and \thetag{1.5} $\Psi(V_t)$ defines
the multiplication by $-1$ on  $\widehat{L}_t/L_t$, i.e.
$V_t$ corresponds  to the element
$\xi_t=-1$ of $\Xi(t)$ (see \thetag{1.4}).
Therefore
$$
-\Psi({V}_t)\in \widehat{\hbox{O}}(L_t).
$$
Elements $M$ and $-M\in \operatorname{O}(L_t)$ define the same
transformation of the domain $\bh_t^+$. Thus we have
\proclaim{Corollary 1.3}Let $t$ be square free. The groups
$$
\Psi(\gt\cup\gt V_t)\qquad\text{and}\qquad
\operatorname{O }^*(L_t)=\widehat{\hbox{\rm O}}(L_t)\cap
\operatorname{O }^+_\br(L_t)
$$
coincide, if  we consider them  as  groups of analytic
transformations of  $\ot $.
\endproclaim
\proclaim{Proposition 1.4}The quotient
$$
(\gt\cup\gt V_t)\setminus \bh_2
$$
is isomorphic to the moduli space of polarized $K3$ surfaces with
a polarization of type $<2t>\oplus\, 2E_8(-1)$.
\endproclaim
\demo{Proof}
It is known that a moduli space of polarized $K3$ surfaces is a quotient
of a $19$-dimensional homogeneous domain of type IV
by an arithmetic group. In the proposition we consider
polarized $K3$ surfaces with a condition on the Picard group
or equivalently on the  lattice of its trancendental cycles.
To formulate these conditions we need some definitions (see
\cite{N2}, \cite{D}).

Let $X$ be a $K3$ surface.
Let us take a sublattice $D_t=<2t>\oplus \,E_8(-1)\oplus E_8(-1)$
of the lattice
$$
L_{K3}=
U\oplus U\oplus U\oplus E_8(-1)\oplus E_8(-1)
\cong H^2(X,\Bbb Z),
$$
where $<2t>$ ($t\in \bz$) denotes the  one-dimensional lattice
generated by a vector $l$
such that $l^2=2t$,
$U$ is the hyperbolic plane with quadratic form
$\left(\smallmatrix 0&1\\1&0\endsmallmatrix\right)$
and $E_8(-1)$ is the even  unimodular lattice of
dimension 8 with the negative definite quadratic form.
We note that
$\hbox{sign\,}(D_t)=(1,16)$ and
$D_{t}^\perp\cong U\oplus U\oplus <-2t>=L'_t=L_t(-1)$.
(Notation $L_t(-1)$ means that we multiply the quadratic form $S_t$
on the lattice $L_t$ by $-1$.)
We recall (see \cite{N1}) that the orthogonal  group
$$
\operatorname{O}(D_t,\,L_{K3})=
\{\,g:L_{K3}\to L_{K3}\ |\ g\,|_{D_t}\equiv \hbox{id}\}
$$
is isomorphic to the group
$$
\widehat{\hbox{O}}(L'_t)=\{\,g:\,L'_t\to L'_t\ |\ \forall\  l\in
\widehat{L}'_t\ \  gl- l \in  L'_t\,\}.
$$

A marked $D_t$-polarized $K3$-surface is defined by the following
datum (see  \cite{N2} and  \cite{D} for more details):
a surface $X$,  a fundamental domain $C(M)^{+}$ of a group
generated by some $2$-reflections of the lattice $D_t$ acting on
a connected  component $V(D_t)^+$ of the cone
$V(D_t)=\{v\in D_t\otimes \br\ |\ (v,v)>0\}$
and  an isomorphism of the lattices
$\phi:\,H^2(X,\bz)\to L_{K3}$,
such that
$\phi^{-1}(D_t)\subset\hbox{Pic\,}(X)$,
$\phi^{-1}(V(D_t)^+)\subset V(X)^+$ and
$\phi^{-1}(C(M)^+)$
contains at least one numerically effective divisor class.
By  $V(X)^+$ one  denotes the connected component of the cone
$$
V(X)=\{v\in H^{1,1}_\br (X)\ |\ (v,\,v)>0\}
$$
containing the cohomology class of a K\"ahler form on $X$.

Let us denote by $\omega_X$  a  holomorphic $2$-form
which generates  $H^{2,0}(X)$.
Its image  under the isometry $\phi$ belongs
 to the following domain in the projective space $\Bbb P^{4}$
$$
\phi(\omega_X)\in \Cal D_m=\{v\in \Bbb P(L'_t\otimes\Bbb C):\ (v, v)=0,
\quad (v,\bar v)>0\}.
$$
This domain is an example of the domains of type IV.
Its connected components are isomorphic to the domain
$\Cal H_t^+$ introduced in   \thetag{1.7}.

The quotient
$$
\Cal M(<2t>\oplus \,2E_8(-1))=
\operatorname{O}^*(L'_t)\setminus \Cal H_t^{+}
$$
is the  moduli space of isomorphism classes of
 $<2t>\oplus\, 2E_8(-1)$-polarized $K3$ surfaces.
In accordance with Corollary 1.3
$$
\Cal M(<2t>\oplus \,2E_8(-1))\cong (\gt\cup \gt V_t) \setminus \bh_2.
$$
\qed\enddemo

We now want to explain the relationship between the variety
${\Cal A}_t^*=\Gamma_t^*\setminus \bh_2$ and the space of
Kummer surfaces associated to $(1,t)$-polarized abelian surfaces.
For an abelian surface $A$ we denote by $X=\km (A)$
its associated Kummer surface.

\proclaim{Theorem 1.5} (i) Let $A$, $A'$ be two $(1,t)$-polarized
abelian surfaces  which define the same point in  ${\Cal A}_t^*$.
Then their Kummer surfaces $X$, $X'$ are isomorphic.

(ii) Assume that the Neron-Severi group of $A$, resp. $A'$ is
generated by the polarization. Then the converse is true:
If $A$ and  $A'$ have isomorphic Kummer surfaces,
then $A$ and  $A'$ define the same point in  ${\Cal A}_t^*$.
\endproclaim
\demo{Proof}We consider the commutative diagram
$$
\CD
\tilde A@>p>>X=\km (A)\\
@V\sigma VV @VVV\\
A@>>>Y=A/\iota
\endCD
$$
where $\tilde A$ is the blow-up of $A$ in the $16$ points of order $2$.
By \cite{BPV, p. 246} the map
$$
\alpha=p_!\circ \sigma^*:H^2(A,\bz)\to H^2(X,\bz)
$$
multiplies the intersection form  by $2$ and hence is in particular a
monomorphism. Let $T\subset H^2(X,\bz)$ be the orthogonal complement of
the sublattice of $H^2(X,\bz)$ generated by the polarization
and the nodal classes. Using \cite{BPV, Corollary VIII.5.6}
it follows that  $T\cong 2U(2)\oplus <-4t>$.
If $A$ is generic, i.e. the Neron-Severi
group is generated by the polarization,
then $T$ is the transcendental lattice of $X$.

(i) Assume that $A$ and $A'$ define the same point in ${\Cal A}_t^*$.
Then there exists an element $V$ in $\Gamma_t^*$ which induces an
isomorphism of the lattices $T$ and $T'$ of $X$ and $X'$ respectively.
By the Torelli theorem for Kummer surfaces it now suffices
to prove
\enddemo
\proclaim{Claim}The isomorphism $\ V:T\to T'$ can be extended to a
Hodge isometry
\newline
$\widetilde V: H^2(X,\bz)\to H^2(X',\bz)$.
\endproclaim
\demo{Proof of Claim}By definition $T$ is a primitive
non-degenerate sublattice of
the unimodular lattice $H^2(X,\bz)$.
Its orthogonal complement $T^\perp$ in  $H^2(X,\bz)$
is an  even indefinite lattice of  signature $(1,16)$.
The discriminant groups $A_T$ and $A_{T^\perp}$ are isomorphic
and have five generators. By \cite{N1, Theorem 1.14.2}
the homomorphism
$\operatorname{O}(T^\perp)\to
\operatorname{O}(A_{T^\perp},q_{T^\perp})$ is surjective.
Using \cite{N1, Corollary 1.5.2} one can, therefore, construct
an extension of the isomorphism $V$.
\enddemo
(ii) An isomorphism $f:X\to X'$ induces an isomorphism of
transcendental lattices. Under the assumption stated
this  defines an  isometry $V:T\to T'$ where
$T\cong T'\cong 2U(2)\oplus <-4t>=L_t(-2)$, which is an element
of the group $\operatorname{O}^+(L_t(-2))$.
Since $\operatorname{O}^+(L_t(-2))=\operatorname{O}^+(L_t)$
it follows from Proposition 1.2 that there is an element
$\overline V\in \Gamma_t^*$ which defines one of the elements
$\pm V$. The element
$\overline V$ identifies the points defined by $A$ and $A'$ resp..
\newline\qed
\remark{Remark}The above Theorem justifies it to consider
$\Cal A_t^*$ as the moduli space of Kummer surfaces
associated to abelian surfaces with a $(1,t)$-polarization.
\endremark
\remark{Remark}In \cite{D, Example 6.5} Dolgachev considered the
space of $M_t$-polarized abelian surfaces where $M_t$ is
the orthogonal complement of $2U(2)\oplus <-4t>$ in $L_{K3}$.
This leads to a subgroup of finite index of $\operatorname{O}^+(L_t)$
and hence a  covering space  of $\Cal A_t^*$.
\endremark

Our next aim is to interprete the involutions  $V_d$ geometrically.
Because of Proposition 1.2 the element  $V_d$ induces a map from
${\Cal A}_t^*$ to itself.

Let $(A,H)$ be a $(1,t)$-polarized abelian surface. The polarization $H$
defines an isogeny
$$
\align
\lambda_H:\quad A   &   \rightarrow     {\widehat A} = \hbox{Pic}^0 A\\
                x   &   \mapsto         T_x^*{\Cal L}\otimes{\Cal L}^{-1}
\endalign
$$
where ${\Cal L}$ is a line bundle representing $H$ and $T_x$
denotes translation by $x$. The map $\lambda_H$ only depends on $H$, not on
the line bundle
${\Cal L}$. There is a (non-canonical) isomorphism
$\hbox{ker\,} \lambda_H\cong
\bz_t \times \bz_t$. For every divisor $d$ of $t$ there is a unique
subgroup $G(d)\subset \hbox{ker\,} \lambda_H$
which is isomorphic to $\bz_d \times
\bz_d$. This subgroup defines a quotient
$$
\lambda_d :\quad A \rightarrow A/G(d)=A'.
$$
If $A$ is given by the period matrix
$$
\Omega = \pmatrix
1   &   0   &   \tau_1   &   \tau_2\\
0   &   t   &   \tau_2   &   \tau_3
\endpmatrix, \qquad \tau =
\pmatrix
\tau_1   &   \tau_2\\
\tau_2   &   \tau_3
\endpmatrix
\in \bh_2
$$
then $A'$ is given by
$$
\Omega ' = \pmatrix
d   &   0   &   d\tau_1   &   \tau_2\\
0   &   t_d &   \tau_2   &   \tau_3/d
\endpmatrix.
$$
The abelian surface $A'$ carries a uniquely determined
polarization $H'$ with
$$
dH=\lambda_d^* (H').
$$
The polarization $H'$ is of type $e\cdot (1,t/e^2)$ where $e=(d,t_d)$.
Altogether this shows that we have a morphism of moduli spaces
$$
\align
\Phi=\Phi(d):\quad {\Cal A}_t   &   \rightarrow {\Cal A}_{t/e^2}\\
(A,H)   &   \mapsto   (A',H').
\endalign
$$
If $d=t$ we obtain as a special case the map
$$
\align
\Phi(t):\quad {\Cal A}_t &     \rightarrow {\Cal A}_t\\
(A,H)   &   \mapsto      ({\widehat A},{\widehat H})
\endalign
$$
which maps an abelian surface to its dual polarized abelian surface.
\proclaim{Proposition 1.6}
Let $d$ be a divisor of $t$ with $(d,t_d)=1$. Then the map
$$
\Phi(d):\quad {\Cal A}_t  \rightarrow {\Cal A}_t
$$
is the map induced by $V_d$.
\endproclaim
\demo{Proof}
For $\tau=
\pmatrix
\tau_1  &  \tau_2\\
\tau_2  &  \tau_3
\endpmatrix
\in \bh_2$ we have the following formula for the action
$$
V_d<\tau>=
\pmatrix
x      &   -1\\
-y t_d &   d
\endpmatrix
\pmatrix
d\tau_1  &   \tau_2\\
\tau_2   &   \tau_3/d
\endpmatrix
\pmatrix
x   &   -y t_d\\
-1  &   d
\endpmatrix.
$$
Now consider the matrix
$$
\pmatrix
1   &   t_d   &   0   &   0\\
y   &    xd   &   0   &   0\\
0   &    0    &   x   &   -y t_d\\
0   &    0    &   -1  &   d
\endpmatrix
\in \operatorname{SL}_4(\bz).
$$
This matrix transforms the symplectic form
$d e_1\wedge e_3 + t_d e_2\wedge e_4$
into $W_t$  (see \thetag{1.1}).
The claim now follows from the equality
$$
\gather
\pmatrix
x   &   -1\\
-y t_d   &   d
\endpmatrix \hskip-1.5pt
\pmatrix
d   &   0   &   d\tau_1    &   \tau_2\\
0   &   t_d   &   \tau_2   &   \tau_3/d
\endpmatrix
\hskip-1.5pt\pmatrix
1   &   t_d   &   0    &   0\\
y   &   xd    &   0    &   0\\
0   &   0     &   x    &   -y t_d\\
0   &   0     &   -1   &   d
\endpmatrix\\
= \left(\matrix
1   &   0\\
0   &   t
\endmatrix\ \ V_d<\tau>\right).
\endgather
$$
\qed\enddemo

\remark{Remark}In view of Theorem 1.5 this shows in particular
that a $(1,t)$-polarized abelian surface and its dual
polarized abelian surface have isormorphic Kummer surfaces.
\endremark

\head
\S\ 2. Nonunirationality of some quotient spaces
\endhead

In accordance
with  Proposition 1.4 a moduli space of special
$K3$ surfaces is a quotient of a moduli space of polarized abelian
surfaces. It gives us a double covering
$$
{\Cal A}_{t}\to ({\Gamma}_{t}\cup {\Gamma}_{t} {V}_{t}) \setminus \bh_2.
$$
The degree of the covering
$
{\Cal A}_{t} \to  {\Cal A}_{t}^{*}={\Gamma}_{t}^{*}\setminus \bh_{2}
$
has order $2^{\nu(t)}$, where $\nu(t)$ is the number of prime divisors
of $t$. Since for square free $t$ the extension
$\hbox{P}\Gamma_{t}^*$ is
the  maximal discrete subgroup of $\operatorname{PSp }_4(\br)$
containing $\hbox{P}\Gamma_t$ the quotient
${\Cal A}_{t}^*$ is the minimal Siegel threefold
associated to the  polarization $(1,t)$.

There are $2^{\nu(t)}-2$
other threefolds between ${\Cal A}_{t}$ and ${\Cal A}_{t}^*$.
Let us take for example  two primes $p\neq q$ and let
$t=pq$. The involutions $V_p$ and $V_q$ give rise to the following moduli
spaces
$$
{\Cal A}^{(p)}_{pq} = {\Cal A}_{pq} /<V_p>,\
{\Cal A}^{(q)}_{pq}= {\Cal A}_{pq} /<V_q>,\
{\Cal A}^*_{pq} = {\Cal A}_{pq} /<V_p,V_q>\,
=\Gamma^*_{pq}\backslash \bh_2
$$
resp. a commutative diagram
$$
\matrix
       &         & {\Cal A}_{pq}&  &\\
       & \swarrow &             & \searrow   &\\
{\Cal A}^{(p)}_{pq} &   &    &  & {\Cal A}^{(q)}_{pq}\\
       & \searrow   &       & \swarrow  & \\
       &         & {\Cal A}^*_{pq} &          &
\endmatrix
$$
where all maps are 2:1.
Using the modular forms constructed in  \cite{G1}
we can obtain  information about the geometrical type
of some of these moduli spaces.

By $J_{k,t}^{cusp}$ we denote the space of Jacobi cusp forms of weight $k$
and index $t$. In \cite{G1} a lifting was constructed which associates to
a Jacobi form $\Phi \in J_{k,t}^{cusp}$  a cusp form
$F_{\Phi}\in {\Cal M}_k(\widehat{\Gamma}_t)$
of weight $k$ with respect to the group
$$
\widehat{\Gamma}_t=
\left\{\pmatrix
*     &  t*   &   *   &   *\\
{*}   &  *    &   *   &   *t^{-1}\\
{*}   &  t*   &   *   &   *\\
t*    &  t*   &  t*   &   *
\endpmatrix
\in \operatorname{Sp }_4
(\bq) \right\},
$$
where all entries $*$ denote integers.

The groups $\Gamma_t$ and $\widehat{\Gamma}_t$ are conjugate.
Indeed if
$C_t=\hbox{diag\,}(1,\,t^{-1},\,1,\,t)$
then
$
\widehat\Gamma_ t=C_t \Gamma_t C_t^{-1}
$.
For
$$
\widehat{V}_t={}^tV_t=
\pmatrix
0   &   \sqrt{t}   &   0   &   0\\
\sqrt{t}^{-1}   &   0   &   0   &   0\\
0   &   0   &   0   &   \sqrt{t}^{-1}\\
0   &   0   &   \sqrt{t}   &   0
\endpmatrix
$$
it was proved in \cite{G1, formula (2.8)} that
$$
F_{\Phi}(Z)=(-1)^k F_{\Phi}|_k\, \widehat{V}_t(Z).
$$
or equivalently
$$
F_{\Phi} (\pmatrix
\tau_1   &   \tau_2\\
\tau_2   &   \tau_3
\endpmatrix) =
F_{\Phi} (\pmatrix
t\tau_3   &   \tau_2\\
\tau_2    &   t^{-1} \tau_1
\endpmatrix).
$$
In this section we describe the behavior of $F_{\Phi}$ with respect
to the group
$$
\widehat\Gamma_t^*=\,<\widehat\Gamma_t, \, \widehat{V}_d\ |\ d||t>,
$$
where $\ \widehat{V}_d=C_t V_d C_t^{-1}$.

Eichler and Zagier \cite{EZ, \S\ 5}, have constructed a decomposition of
the space of Jacobi forms
$$
J_{k,t}^{cusp}=\bigoplus\limits_{\epsilon}\, J_{k,t}^{\epsilon}
$$
where $\epsilon$ runs over all characters of the group
$$
\Xi (t)=\{\,\xi \operatorname{mod\,} 2 t\ |\  \xi^2\equiv 1
\operatorname{mod\,}4t\,\}\cong (\bz/2\bz)^{\nu (t)}
$$
satisfying $\ \epsilon(-1)=(-1)^k$.
For any $d$ with $d || t$ one can define an operator $W_d$
acting on $J_{k,t}$ in the following way \cite{EZ, \S 5}.
For
$$
\Phi(\tau_1,\tau_2)=\sum\Sb n,\,l\in \bz\\ 4nt >l^2 \endSb
f(n,l) \exp (2\pi i (n\tau_1 + l\tau_2))\in J_{k,t}^{cusp}
$$
we put
$$
(\Phi\,|\,W_d)(\tau_1,\tau_2)=\sum\Sb n,\,l\in \bz\\ 4nt > l^2 \endSb
f(n', l') \exp (2\pi i (n \tau_1 + l\tau_2))\in J_{k,t}^{cusp}
$$
where $l', n'$ are determined by
$$
l'\equiv -l \operatorname{\,mod\,} 2d,
\quad l'\equiv l \operatorname{\,mod\,} 2t/d,
\quad 4n' t-{l'}^{2}= 4nt-l^2.
$$
All $W_d$ are involutions. They form a group isomorphic to $\Xi(t)$.
The subspaces $J_{k,t}^{\epsilon}$ are eigenspaces of the operation
$W_d$, namely
$$
J_{k,t}^{\epsilon}=\{\Phi \in J^{cusp}_{k,t}\ |\  \Phi\,|\, W_d
= \epsilon(W_d)\,\Phi\}.
$$
Note that if
$$
\Phi(\tau_1, \tau_2)=\sum_{\mu \operatorname{mod} 2 t} \varphi_{\mu} (\tau_1)
\theta_{t,\mu}(\tau_1,\tau_2)\in J_{k,t}^{\epsilon}
$$
is the standard decomposition of the Jacobi form $\Phi$ with respect to
the theta-functions $\theta_{t,\mu}(\tau_1, \tau_2)$ then for
$\xi_d\in \Xi(t)$
$$
\varphi_{\xi_d\mu}(\tau_1)=\epsilon(\xi_d) \varphi_{\mu} (\tau_1).
$$
\proclaim{Theorem 2.1}
Let $\Phi\in J_{k,t}^{\epsilon}$ be a Jacobi form and $F_{\Phi}\in {\Cal
M}_k (\widehat{\Gamma}_t )$ be its lifting. For any divisor $d$ of $t$ with
$(d, t_d)=1$ the following equality holds
$$
F_{\Phi}|_k \,\widehat{V}_d = \epsilon (\xi_d) F_{\Phi}.
$$
\endproclaim
\demo{Proof}
Let us recall the definition of the lifting $F_{\Phi}$ in terms of the
Fourier expansion \cite{G2}. If
$$
\Phi(\tau_1,\tau_2)=\sum\Sb n,\,l\in \bz\\ 4nt > l^2 \endSb
f(n, l) \exp (2\pi i (n \tau_1 + l\tau_2))\in J_{k,t}^{cusp}
$$
then
$$
F_{\Phi}(Z)=\sum_{N\in \frak{A}_t} b(N)\exp (2\pi i\ tr (NZ))
$$
where summation is taken over all positive definite symmetric matrices
of the following form
$$
N\in {\frak{A}_t}=\left\{  \pmatrix
n     &   l/2\\
l/2   &   mt
\endpmatrix > 0\ |\  n, l, m \in \bz \right\}
$$
and
$$
b( \pmatrix
n     &   l/2\\
l/2   &   mt
\endpmatrix)= \sum_{a\,|\,(n, l, m)} a^{k-1}\
f \left(\frac {n m}{ a^2}, \frac l a \right).
$$
The action of $\widehat{V}_d$ on $F_{\Phi}$ is given by
$$
(F_{\Phi} |_k\, \widehat{V}_d)(Z)
=F_{\Phi}(d^{-1} A_d\ Z\ {} ^tA_d)\quad \hbox{ where }
A_d= \pmatrix
dx   &   -t\\
-y   &   d
\endpmatrix
$$
or
$$
\align
(F_{\Phi} |_k \,\widehat{V}_d)(Z)&  =  \sum_{N\in\frak{A}_t} b(N)
\exp(2\pi i\ tr (d^{-1}\ {} ^tA_d N A_dZ))\\
{}& =  \sum_{N\in\frak{A}_t} b(d^{-1}\ {} ^t{\tilde A}_d N
{\tilde A}_d) \exp (2 \pi i\ tr (NZ))
\endalign
$$
where
${\tilde A}_d = d A_d^{-1} =
\pmatrix
d   &   t\\
y   &   dx
\endpmatrix$.
Let
$
{\widetilde N}= d^{-1}\ {} ^t{\tilde A}_d N {\tilde A}_d =
\pmatrix
{\tilde n}   &   {\tilde l}/2\\
{\tilde l}/2 &   {\widetilde m} t
\endpmatrix
$.
Clearly $\hbox{det\,} N =\hbox{det\,} {\widetilde N}$.
It is easy to see that the elements
$n$, $l$, $m$ and ${\tilde n}$, ${\tilde l}$, ${\widetilde m}$
have the same set of
common divisors. Moreover
$$
{\tilde l} = l (y t_d + dx) + 2 (nt + xymt)\quad \hbox{and}
\quad{\tilde l}\equiv
\cases -l &\operatorname{\,mod\,} 2d\\
       \hphantom{-}l &\operatorname{\,mod\,} 2t_d
\endcases
$$
Hence, by the definition of $J_{k,t}^{\epsilon}$  we have
$$
b(d^{-1}\ {} ^t{\tilde A}_d N {\tilde A}_d)=\sum_{a\,|\,(n, l, m)}
a^{k-1}  f\left( \frac {{\tilde n}{\widetilde  m}}{a^2},\,
\frac {{\tilde  l}}{a}\right) =
\epsilon (\xi_d) b(N)
$$
which proves the theorem.
\newline\qed\enddemo

This result can be used to gain some information on
the Kodaira dimension of moduli spaces.
Whenever we speak of the {\it Kodaira dimension} of some moduli space
${\Cal A}$ we mean the Kodaira dimension of a desingularization of a
projective compactification of ${\Cal A}$.

\proclaim{Corollary 2.2}
Let $p\neq q$ be primes $\ge 5$. Then the Kodaira dimension of at least
one of the spaces ${\Cal A}^{(p)}_{pq}$  or ${\Cal A}^{(q)}_{pq}$ is
$\ge 0$.
\endproclaim
\demo{Proof}
The Eichler-Zagier decomposition gives a decomposition
$$
J_{3, pq} = J_{3, pq}^{(+,-)} \oplus J _{3, pq}^{(-,+)}.
$$
If $p,q \ge 5$ then there exists a cusp form in $J _{3, pq}$ and hence
in $J_{3, pq}^{(+,-)}$ or $J_{3, pq}^{(-,+)}$. By \cite{G1},
\cite{G2} this
can be lifted to a weight 3 cusp form with respect to
$<\Gamma_t,\, V_p>$, resp. $<\Gamma_t,\, V_q>$.
By Freitag's extension theorem this
defines a differential form on any desingularization of a projective
compactification of ${\Cal A}^{(p)}_{pq}$, resp. ${\Cal A}^{(q)}_{pq}$.
This gives the result.
\newline\qed\enddemo

For any integer $t$
let us take a character of the group $\Xi(t)$ isomorphic to
the orthogonal group of the  discriminant group  of $L_t$
(see \thetag{1.4})
$$
\epsilon:\ \Xi(t)\to \{\pm 1\}.
$$
We define a set $U(\epsilon)=\{\,V_d\in \gts |\ \epsilon(\xi_d)=1\,\}$
and a subgroup
$$
\gt^{\epsilon}=\,<\gt,\ \xi_d\ |\ \xi_d\in U(\epsilon)>\,\subset \gts
$$
of $\gts$.
Theorem 2.1 and the method of the proof of Corollary 2.2 gives us
the next result
\proclaim{Corollary 2.3}If $\hbox{\rm dim\,}(J_{3,t}^{\epsilon})>0$,
then the Kodaira dimension of the   quotient space
$$
\Cal A_t^{\epsilon}= \gt^{\epsilon}\,\backslash \bh_2
$$
is nonnegative.
\endproclaim
\remark{Remark}Corollary 2.3 gives information about
$\Cal A_t^{\epsilon}$ only if $V_t\notin U(\epsilon)$.
If $\epsilon(V_t)=\epsilon(-1)=1$, then
$\hbox{dim\,}(J_{3,t}^{\epsilon})=0$.
\endremark
\proclaim{Corollary 2.4}Let $t\ge 21$ ($t\ne 30,\, 36$)
and let its number of prime divisors
$\nu(t)\ge 2$.
Then there exists a finite quotient of
$\Cal A_t$ of degree $2^{\nu(t)-1}$ which is not unirational.
\endproclaim
\demo{Proof}For any integer $t$ from the corollary the dimension
of $J_{3,t}^{cusp}$ is positive. Thus there is a character
$\epsilon$ of $\Xi(t)$ such that $\epsilon(\xi_t)=-1$
and $\hbox{dim}\,J_{3,t}^{\epsilon}>0$.
$\Xi(t)\cong (\bz/2\bz)^{\nu(t)}$ therefore
$[\Gamma_t^{\epsilon}:\gt]=2^{\nu(t)-1}$.
\newline\qed\enddemo
Using dimension formulae for the spaces $J^{\epsilon}_{k,t}$
one can obtain more precise results.
It is easy to get an exact dimension formula  using
the  trace formula of the operator $W_d$ on the space
$J_{k,t}^{cusp}$ given in \cite{SZ}.  By definition of
$W_d$ we have
$$
\hbox{tr\,}(W_d,\  J_{3,t}^{cusp})=
\sum_{\epsilon}\,\epsilon(\xi_d)\,
\hbox{dim\,}(J_{3,t}^{\epsilon}),
$$
where the sum is taken over all characters of $\Xi(t)$.
Therefore
$$
\hbox{dim\,}(J_{3,t}^{\epsilon})=
2^{-\nu(t)}\sum_{d||t}\epsilon(\xi_d) \,
\hbox{tr\,}(W_d,\ J_{3,t}^{cusp}).
$$
For weight $3$ the trace formula  of  $W_d$ on $J_{3,t}^{cusp}$
(we recall that $d|t$ and $(d,t_d)=1$, where $t_d=\dsize\frac t{d}$\,)
proved in  \cite{SZ, Theorem 1} can be  reduced to the following
expression
$$
\multline
\hbox{tr\,}(W_d,\ J_{3,t}^{cusp})=
\frac 1{4}\sum_{e|d}H_{t_d}(-4e)-\frac 1{4}\sum_{e'|t_d}H_{d}(-4e')+
\frac 3{2}(H_d(0)-H_{t_d}(0))
\\
+\frac 1{2}\bigr(\delta_2(t_d)H_d(-4)-\delta_2(d)H_{t_d}(-4)\bigl)+
\bigr(\delta_3(t_d)H_d(-3)-\delta_3(d)H_{t_d}(-3)\bigl)\\
+\frac 1{4}\biggr((Q(t_d),2)Q(d) -(Q(d),2)Q(t_d)\biggl).
\endmultline
$$
We denote by $Q(n)$ the greatest integer whose square divides $n$;
$\delta_a(b)=1$ or $0$ if $a|b$ or $a\not| \ b$
and $H_n(\Delta)$ is a generalization of the Hurwitz-Kronecker
class number, i.e $H_1(0)=-\frac 1{12}$
and $H_1(\Delta)$ for $\Delta<0$ is the number of equivalence classes
with respect to $\operatorname{SL}_2(\bz)$ of integral, positive definite,
binary quadratic forms of discriminant $\Delta$,
counting forms equivalent to a multiple of $x^2+y^2$
(resp. $x^2+xy+y^2$) with multiplicity $\frac 1{2}$ (resp. $\frac 1{3}$).
For $n\ge 2$ with  $(n,\Delta)=a^2b$ and square free $b$
$$
H_n(\Delta)=\cases
a^2b\biggl(\dsize\frac{\Delta/a^2b^2}{n/a^2b}\biggr)
H_1(\Delta/a^2b^2)
\quad&\text{if } a^2b^2|\Delta\\
0&\text{otherwise,}
\endcases
$$
where $\bigl(\frac{\cdot}{\cdot}\bigr)$ is the generalized
Kronecker symbol.

We note that the trace formula has the simplest form for  square free
$t$ coprime to $6$:
$$
\hbox{tr\,}(W_d,\ J_{3,t}^{cusp})=
\frac 1{4}\biggl(\sum_{e|d}\biggl(\frac{-4e}{t_d}\biggr)H_{1}(-4e)-
\sum_{e'|t_d}\biggl(\frac{-4e'}{d}\biggr)H_{1}(-4e')\biggr)+
\frac{t_d-d}{8}.
$$
\example{Example 2.5}The calculation gives us only
thirteen different threefolds
of type $\Cal A_t^{\epsilon}$ with $t$ having only two prime divisors
($t=p^a q^b$),
whose geometric genus could be equal to zero.
They  correspond  to the following {\it trivial\,} subspaces of
$J_{3,t}^{cusp}$ of type $J_{3,p^a q^b}^{\ \,\,-\,+}$
(this notation means that
$J_{3,p^a q^b}^{\ \,\,-\,+}=J_{3,p^a q^b}^{\epsilon}$ with
$\epsilon(\xi_{p^a})=-1$  and $\epsilon(\xi_{q^b})=1$):
$$
\gather
J_{3,2\cdot 11}^{\ \,-+},\ \
J_{3,2\cdot 13}^{\ \,-+},\ \ J_{3,2\cdot 17}^{\ \,-+},\
\ J_{3,2\cdot 19}^{\ \,-+},\ \ J_{3,2\cdot 25}^{\ \,-+},\ \
J_{3,2\cdot 27}^{\ \,-+},\\
J_{3, 3\cdot 7}^{\ \,-+},\ \  J_{3, 3\cdot 13}^{\ \,-+},\ \
J_{3, 3\cdot 16}^{\ \,-+},\ \
J_{3, 5\cdot 7}^{\ \,-+},\ \  J_{3, 5\cdot 8}^{\ \,-+},\ \
J_{3, 7\cdot 4}^{\ \,-+},\ \ J_{3, 7\cdot 8}^{\ \,-+}.
\endgather
$$
We may add to this list twenty  threefolds $\Cal A_t$
with
$t= 1,\,\dots\,,12$, $14$, $15$, $16$, $18$, $20$, $24$,
$30$, $36$ (see \cite{G2})
whose geometric genus   could be zero. According to classical
results and new results of M. Gross and S. Popescu (see \cite{GP})
it really is  for  $t=1,\dots\,,12$, $14$, $16$, $18$, $20$.
\endexample
\example{Example 2.6}We have the  following quotients of order
$4$ and $8$ which are not unirational. One has
$$
\hbox{dim}\,J_{3,42}^{cusp}=
\hbox{dim}\,J_{3,2\cdot 3\cdot 7}^{\ \,++-}=1.
$$
Thus
for
$\Cal A_{42}^{(2,3)}=<V_2, V_3>\setminus \,\Cal A_{42}$
we obtain $h^{3,0}(\Cal A_{42}^{(2,3)})\ge 1$.

The geometric genus of all four threefolds of  type
$<V_a, V_b, V_c>\setminus\, \Cal A_{210}$, where
$a,\ b,\ c\in \{ 2,\,3,\,5,\,7 \}$ is positive.
\endexample
\vskip0.5truecm

\head
\S\  3. The  Humbert surfaces and the ramification locus
\endhead

If one wants to determine the Kodaira dimension of the variety
${\Cal A}^*_t=\gts\backslash \bh_2$ it is important to
know the ramification locus of the covering map ${\Cal A}_t\rightarrow
{\Cal A}^*_t$, i.e. the locus where the stabilizer of the finite group
$\gts/ \gt$ is not trivial.
Unfortunately this turns out to be a difficult question.
Here we shall give a partial answer, i.e. we shall determine
the divisorial part of the ramification locus which is a union of
a finite numbers of Humbert surfaces. We shall restrict ourselves
to $t$ square free.

First  we collect some known facts about  divisors
on the homogeneous domain $\ot$.
For any  $v\in L_t\otimes\br$ we set
$$
\Cal H_v=\{Z\in \ot \ |\ v\cdot Z=0\},
$$
where $v\cdot u=(v, u)_t$ is the bilinear product
corresponding to the quadratic form $S_t$ (see \thetag{1.2}).

\proclaim{Lemma 3.1}

{\rm 1}. $\Cal H_{g v}=g^{-1}\Cal H_v$ for any
$g\in \hbox{\rm O}_\br^+(L_t)$.

{\rm 2.} Let $v\not= 0$, then $\Cal H_v\not= \varnothing$
if and only if $v^2>0$.

{\rm 3}. $\Cal H_{v}\cap \Cal H_{u} \not=\varnothing$
if and only if the matrix
$\pmatrix v^2&v\cdot u\\v\cdot u&u^2\endpmatrix$
is positive define.
\endproclaim
\demo{Proof} 1. The first property is trivial.

2. The orthogonal group $\operatorname{O }^+_\br(L_t)$
acts transitively on $\ot$.
Thus any $Z\in \ot$ can be reduced to
$Z_{\bold i}={}^t(1,i,0,i,1)$ in the   coordinates
$(z_i)$ from \thetag{1.7}.
If $v=(a,b,c,d,e)$ and $v\cdot Z_{\bold i}=0$,
then $a=-e$ and  $b=-d$.
Thus $v^2=2a^2+2b^2+2tc^2>0$.

Let $L_t\otimes \br=\br v\oplus V$ with $v^2>0$. One has
$\hbox{sign\,}(V)=(2,2)$.
The group
$\operatorname{SO }^+(V)\cong \operatorname{SO }_\br^+(2,2)$
is  locally isomorphic to
$\operatorname{SL}_2(\br)\times \operatorname{SL}_2(\br)$.
Thus
$$
\Cal H_v\cong \bp\bh_{V}^+\cong \bh_1\oplus \bh_1,
$$
where $\bh_1$ is the usual upper half-plane.
This proves the second statement.

3. Let us suppose that
$\Cal H_{v}\cap \Cal H_{u} \not=\varnothing$.
It  follows from  $(xu-v)^2>0$ that the matrix in 3 is
positive definite.

If the  symmetric bilinear form on the   plane
$P=\br v\oplus \br u\ $ is
positive definite,
then $\hbox{sign\,}(P^{\perp})=(1,2)$.
The group
$\operatorname{SO }^+(P^{\perp})\cong
\operatorname{SO }^+(1,2)_{\br}$
is  locally isomorphic to $\operatorname{SL}_2(\br)$
and
$$
\Cal H_{v}\cap \Cal H_{u}\cong \bp\bh^+_{P^{\perp}}\cong \bh_1.
$$
\qed\enddemo

\remark{Remark}For $l\in L_t$ ($l^2>0$)
the group $\widehat{\hbox{SO}}^+\hskip -2pt(l^{\perp})$
is isomorphic  to a subgroup
of $\operatorname{SL}_2(\bz)\times \operatorname{SL}_2(\bz)$ or
to   a subgroup of a Hilbert modular group.
\endremark
\definition{Definition}Let $\ell\in \widehat L_t$ be a vector
in  the dual lattice.
The {\it Humbert surface} $H_\ell$ is defined by
$$
H_\ell =\pi\,( \bigcup\Sb
g\in \widehat{SO}^+\hskip-1pt (L_t)\endSb \Cal H_{g \ell}\,),
$$
where $\pi:\, \ot\to \widehat{\hbox{SO}}^+\hskip -2pt(L_t)\setminus\ot$
is the natural projection.
\enddefinition

$\Cal H_\ell$ depends only on the one dimensional lattice
$\Bbb Z \ell$, thus we can restrict ourselves to  primitive
vectors $\ell \in \widehat{L}_t$. The primitivity means that
$\ell/d\not\in \widehat{L}_t$ for any interger  $d>1$.
The first statement  of Lemma 3.1 says that there is a one to one
correspondence between the
$\widehat{\hbox{SO}}^+\hskip -2pt(L_t)$-orbits
of primitive vectors $\ell\in \widehat{L}_t$
with positive norm and  the Humbert surfaces.

It is well known that for any even integral lattice $L$
with two hyperbolic planes
(in particular for $L_t$) the
$\widehat{\hbox{SO}}(L)$-orbit  of any  $l\in L$ depends only on
the norm of $l$  and  its canonical  image
$l^*:=l/\hbox{div\,}(l)$ in
the discriminant group $\widehat L/L$,
where the  {\it  divisor}
$\hbox{div\,}(l)\in \bn$ of $l$ is the positive generator
of the ideal
$\{(x,l)_L\ |\  x\in L\}$. As a corollary we have

\proclaim{Lemma 3.2}Let $\ell_1$, $\ell_2\in \widehat L_t$ be two
primitive  vectors
with the same image in the discriminant group
(i.e. $\ell_1-\ell_2\in L_t$). If $\ell_1^2=\ell_2^2$, then
$H_{\ell_1}=H_{\ell_2}$.
\endproclaim
\demo{Proof}If $\ell_1-\ell_2\in L_t$, then
$\hbox{div\,}(2t\ell_1)=\hbox{div\,}(2t\ell_2)$
and
$\widehat{\hbox{SO}}^+\hskip -2pt(L_t)\,\ell_1
=\widehat{\hbox{SO}}^+\hskip -2pt(L_t)\,\ell_2$.
\newline\qed\enddemo

\definition{Definition}Let $\ell$ be a primitive vector of
$\widehat{L}_t$.
The integer  $\Delta(\ell)=2t\ell^2$ is  called
the {\it discriminant} of $H_\ell$.
\enddefinition

{}From the isomorphism  $\widehat{L}_t/L_t\cong \bz/ 2t\bz$
one gets
\proclaim{Corollary 3.3}The number of  surfaces $H_\ell$
with fixed discriminant $\Delta=2t\ell^2$,
which are not  $\gt$-equivalent,
is equal to the number of solutions
$$
\#\,\{b\operatorname{mod } 2t\,
|\ b^2\equiv \Delta\operatorname{mod } 4t\}.
$$
\endproclaim

The standard definition of the Humbert surfaces
(see \cite{vdG}, \cite{F})
is given in terms of the moduli space of abelian surfaces
with polarization $(1,t)$.
Let us compare  both definitions.

According to \thetag{1.8} and Proposition 1.2
we may rewrite the definition of $\Cal H_\ell$
with  $\ell=(e, a, -\frac{b}{2t}, c, f)\in \widehat{L}_t$
($(e,a,b,c,f)=1$) in  coordinates $(\tau_i)$ of $\bh_2$:
$$
\Cal H'_x=\psi_t^{-1}(\Cal H_\ell)=\{\,
\pmatrix \tau_1&\tau_2\\ \tau_2&\tau_3\endpmatrix\in \bh_2\ |\
(\tau_2^2-\tau_1\tau_3)f+c\tau_3+b\tau_2+ta\tau_1+te=0\,\},
$$
where $x=(te, ta, b, c, f)$.
The number
$$
2t\ell^2=b^2-4f(te)-4c(ta)=\Delta(\Cal H'_x).
$$
is by definition the discriminant of $\Cal H'_x$.
Let us introduce a lattice
$$
N_t=\{(e,a,b,c,f)\in \bz^5\,|\   e,\,a
\equiv 0\operatorname{mod } t\}.
$$
In accordance with
Proposition 1.2 and Lemma 3.1 we have the following decomposition of
the usual (in sense of \cite{vdG}, \cite{F}) Humbert surface
$H_\Delta\subset \Cal A_t=\gt\setminus \bh_2$:
$$
H_\Delta=\pi_t\,\bigl(\bigcup
\Sb x\in N_t,\ primitive\\ \vspace{1\jot}\Delta(x)=\Delta\endSb
\Cal H'_x\,\bigr)\cong
\pi\,\bigl( \bigcup
\Sb 2t\ell^2=\Delta \endSb
\ \bigcup \Sb
g\in \widehat{SO}^+\hskip -1pt(L_t)\endSb
\Cal H_{g \ell}\,\bigr),
$$
where one takes the summation over representatives $\ell$
from the distinct  orbits
and $\pi_t$ is the natural projection $\pi_t:\bh_2\to \Cal A_t$.
Thus the surface  $\Cal H_\ell$  defined above corresponds to
an irreducible component of the  surface $H_\Delta$.
 Corollary 3.3 tells us  that
the number of  the irreducible components of the $H_\Delta$ is equal to
$$
\#\,\{b\operatorname{mod } 2t\,|\ b^2\equiv
\Delta\operatorname{mod } 4t\}.
$$
This gives a new proof of Theorem 2.4 in \cite{vdG} (see p. 212).

In \S\ 1 we  fixed a basis of the lattice $L_t$ such that
$$
L_t=U(-1)\oplus U(-1)\oplus <2t>,
$$
where $U(-1)$ is the integral hyperbolic plane with the quadratic form
$\left(\smallmatrix 0&-1\\-1&0\endsmallmatrix\right)$ and $<2t>$ is
the  one dimensional
$\bz$-lattice
with  even quadratic form $2t$.
By
$L^{(3)}_t$
we denote the orthogonal component in $L_t$
of the first hyperbolic plane
$$
L^{(3)}_t=\bigl(e_2\wedge e_3,\, e_1\wedge e_3 -
te_2\wedge e_4,\,e_4\wedge e_1\bigr)\bz^3\subset L_t.
$$
It is easy to see that in  any orbit  $\widehat{\hbox{SO}}(L_t)\ell$
there is a vector from $L^{(3)}_t$ and
\newline
$\widehat{\hbox{SO}}^+\hskip -2pt(L_t)\ell=\widehat{\hbox{SO}}(L_t)\ell$.
Thus
any Humbert surfaces can be given in the form
$$
\align
\Cal H_\ell&=\{a z_1+b z_2+c z_3=0\}\subset \ot,
\qquad \Delta(\Cal H_\ell)=2t \ell^2\\
\intertext{or}
\Cal H'_{x}&=\{ta\tau_1+b\tau_2+c\tau_3=0\}\subset \bh_2,\qquad
\Delta(\Cal H'_x)=b^2-4tac=2t \ell^2
\endalign
$$
where $\ell={}^t(0,a,-\frac{b}{2t},c,0)\in \widehat{L}_t$  and
$x=(ta,b,c)\in N_t$.

\medskip

For any $d||t$ we define the following subgroup of $\gts$
and the  corresponding  quotient space
of the moduli space $\Cal A_t$
$$
\gt^{(d)}={\gt\cup \gt V_{d}},\qquad \qquad
\Cal A_t^{(d)}=\gt^{(d)}\setminus \bh_2.
$$
The ramification locus of the map
${\Cal A}_t\rightarrow {\Cal A}^{(d)}_t$
can consist of components of different dimension. In the next theorem
we describe  its  {\it divisorial part} $D_t^{(d)}$.

\proclaim{Theorem 3.4}Let $t$ be square free, $d>1$ and
$t_d=\dsize\frac{t}{d}$. Then
$$
D_t^{(d)}=
\cases H_{4d}\cup H_{d}
&\hbox{\rm if \ } \ \biggl(\dsize\frac {d}{4t_d}\biggr)=1\\
H_{4d} & \hbox{\rm if \ } \ \biggl(\dsize\frac {d}{4t_d}\biggr)\ne 1
\ \ \text{\rm  and }\  \biggl(\dsize\frac {d}{t_d}\biggr)=1,
\endcases
$$
where $ \biggl(\dsize\frac {a}{b}\biggr)$
is the generalized Kronecker symbol of
the quadratic residue.
\endproclaim
\remark{Remarks}1. For $d=t$
$$
D_t^{(t)}
=\cases H_{4t}\cup H_{t}&
\hbox{if \ } t\equiv 1\operatorname{mod }4\\
 H_{4t} &
\hbox{otherwise}.
\endcases
$$
In particular $D_t^{(t)}$ is irreducible if $t\equiv 2$ or
$3\operatorname{mod } 4$ (see Corollary 3.3).

2. For $d=1$  Theorem 3.4 is still true if we denote
by  $D_t^{(1)}$  the divisorial part of the branch locus
of the covering $\bh_2\to \Cal A_t$. We note that
$D_t^{(1)}$  was  found in \cite{Br} by another method.
\endremark

\proclaim{Corollary 3.5}Let $t$ be square free. The divisorial part
$D^*_t$ of the ramification locus  of the map
${\Cal A}_t\rightarrow \Cal A_t^*$,
where
$\Cal A_t^* =\gts\setminus \bh_2$
is the ``minimal'' Siegel modular threefold
corresponding  to polarization of type $(1,t)$,
is the  union  of the following Humbert surfaces
$$
D^*_t=\bigcup_{d|t}\,\bigl( \varepsilon_1(d) H_{4d} \cup
\varepsilon_2(d) H_{d}\bigr),
$$
where $\varepsilon_1(d)=1$ if $\biggl(\dsize\frac d{t_d}\biggr)=1$,
$\varepsilon_2(d)=1$ if $d$ is odd and
$\biggl(\dsize\frac d{4t_d}\biggr)=1$
and they  equal $0$ in all other cases. Moreover none of
the above  Humbert surfaces are  $\gt$-equivalent.
\endproclaim
\demo{Proof of Corollary}We have to prove only the last statement.
If $\ell_1$ and $\ell_2\in\widehat{L}_t$ are two
primitive vectors with norms $\ell_1^2=2/t_{d_1}$,
$\ell_2^2=1/(2t_{d_2})$, then $\ell_1^2\not=\ell_2^2$, since
$t$ is square free.
\newline\qed\enddemo

We break up the proof of Theorem 3.4 into several lemmas.

Let us consider a reflection with respect to a vector
$v\in L_t\otimes \br$:
$$
\sigma_v (x)=x-\frac{2(x,v)}{(v,v)}\,v.
$$
It is known that $\sigma_v\in \operatorname{O}_\br^+(L_t)$ if and only if
$v^2>0$. (This follows from the definition of the real spin norm.)
If $\sigma_v\in \operatorname{O}_\br^+(L_t)$,
then the set  $\hbox{Fix\,}_{\ot}(\sigma_v)$
of  fix points of  $\sigma_v$ on $\bp\bh_t^+$ is a complex surface
$\Cal H_v$. The opposite statement is also true.

\proclaim{Lemma 3.6}Let us suppose that the set of  fix points of
$\sigma\in \operatorname{SO }^+_\br(L_t)$ on $\ot$ is a complex surface.
Then $-\sigma$ is a reflection with respect to a vector
$v\in L_t\otimes \br$.
\endproclaim
\demo{Proof}Over $\br$ one can reduce the quadratic form $S_t$
to  $S=\hbox{diag\,}(E_3,-E_2)$.
The maximal compact subgroup $K_\br$ of the orthogonal group
$\operatorname{SO }_\br^+(S)$
is isomorphic to $\operatorname{SO }(3)\times \operatorname{SO }(2)$
consisting of all elements  which
fix the point $Z_{\bold i}={}^t(0,0,0,i,1)\in \bp\bh_t^+$.
Since the group $\operatorname{SO }^+_\br(S)$ acts
transitively on the homogeneous domain
we can  suppose  that   $\sigma=\hbox{diag\,}(A,B)\in K$ where
$A\in \operatorname{SO }(3)$ and $B\in \operatorname{SO }(2)$.
If $B\not=\pm E_2$, then $B$  has only
one fix point ${\bold i}={}^t(i,1)$ on the projective line.
If $\sigma=\hbox{diag\,}(A,B)$
has at least three fixed points, then   $B$ has an
eigenvalue $\lambda$ of order two. $A$ and $B$ are orthogonal,
thus all eigenvalues of $\sigma$ are equal to $\pm 1$.

There are two possibilities for the  set of eigenvalues of $\sigma$
$$
\{\,\lambda(\sigma)\,\}=\{\,1, -1, -1, -1, -1\}\quad\text{or}
\quad \{\,1, 1, 1, -1, -1\}.
$$
In the first case $-\sigma$ is  a reflection. In the second
case $\sigma$ can be written as a product of two reflections
$\sigma_v \sigma_u$
with orthogonal $u$ and $v$. Thus
$\hbox{Fix\,}_{\bp\bh_t^+}(\sigma)=\Cal H_u\cap \Cal H_v$
and  we have proved the lemma
for non-trivial $B$.

If $B=\pm E_2$, then the same arguments show that  $\sigma$ is conjugate
to
$$
D=\pmatrix B_1&0&0\\0&\pm 1&0\\0&0&\pm E_2\endpmatrix \
\qquad B_1\in \operatorname{SO }(2)             \tag{3.1}
$$
if $\sigma$ has at least two fixed points. If $B_1\not=\pm E_2$,
then $\hbox{Fix\,}_{\bp\bh_t^+}(D)$ is a subset of
$\Cal H_x\cap \Cal H_y$,
where $x$ and $y$ form an orthogonal basis of the plane
of rotation of   $B_1\in \operatorname{SO }(2)$.
\newline\qed\enddemo

\proclaim{Lemma 3.7}There is a one to one correspondence between
the irreducible components $H$ of the divisorial part $D_t^{(d)}$
and the surfaces
$H_\ell$ defined by  reflections $\sigma_\ell\in \gt V_d$.
\endproclaim
\demo{Proof}By $\Cal H$ we denote an irreducible  surface in $\bh_2$
whose image is $H$.
Let us suppose that  $\Cal H=\hbox{Fix\,}_{\bh_2}(G)$ with $G\in\gt V_d$.
In accordance with Proposition 1.2 and Lemma 3.6
$\Psi(G)=\sigma_\ell$ is a reflection.
Moreover $\psi_t(\Cal H)=\Cal H_\ell\subset \ot$ and  $\sigma_\ell$
induces  multiplication by
$\xi_d$ on the discrimnant group $A_{t}$.

The reflection
$$
\sigma_\ell(x)=x-\frac{2(x,\ell)}{(\ell,\ell)}\,\ell
$$
depends only on the line $<\ell>$ defined by $\ell\in L_t$.
It follows from the definition  that $\sigma_\ell$
keeps the lattice $L_t$ invariant if and only if
$\ell^2\,|\,2D$ where  $D=\hbox{div\,}(\ell)$
(see the definition before Lemma 3.2).
The surface $H$ depends only on the class
$\{\,\gamma G\gamma^{-1}\,|\,\gamma\in \gt\}$ and
$$
\Psi(\gamma G \gamma ^{-1})= \beta \sigma_\ell \beta^{-1}=
\sigma_{\beta \ell}
\qquad \bigr(\gamma \in \gt,\ \beta=\Psi(\gamma )
\in \widehat{\hbox{SO}}^+\hskip -2pt(L_t)\bigl).
$$
\qed\enddemo
Therefore in order to find all Humbert surfaces in the divisorial part
we have to classify the $\widehat{\hbox{SO}}(L_t)$-orbits of vectors
$\ell\in L_t$ with the additional condition
$\ell^2|\,2\hbox{div\,} (\ell)$.

\proclaim{Lemma 3.8}Let $t$ be an arbitrary positive integer and $d||t$.
There is a one to one correspondence between the
$\widehat{\hbox{\rm SO}}(L_t)$-conjugacy classes of
reflections $\sigma_\ell$ in  the coset  $(-\gt V_d)$
and the   orbits of  the primitive vectors in $L_t$, which satisfy
the following conditions:
$$
\ell^2=2d \quad\text{ and}\quad
\cases
\hbox{\rm div\,}(\ell)=2t_d&\hbox{ if \ }
\biggl(\dsize\frac {d}{4t_d}\biggr)=1
\\
\hbox{\rm div\,}(\ell)=t_d & \hbox{ if \ }
\biggl(\dsize\frac {d}{t_d}\biggr)=1.
\endcases
$$
\endproclaim
\demo{Proof}
One can suppose that  $\ell={}^t(0,a,b,c,0)\in L_t$
and $(a,b,c)=1$. For such $\ell$ the matrix of $\sigma_\ell$
has the following form
$$
\sigma_\ell=\pmatrix
1&0&0&0&0\\
0&1+\frac{2ca}{\ell^2}&-\frac{4tba}{\ell^2}&\frac{2a^2}{\ell^2}&0\\
0&\frac{2cb}{\ell^2}&1-\frac{4tb^2}{\ell^2}&\frac{2ab}{\ell^2}&0\\
0&\frac{2c^2}{\ell^2}&-\frac{4tbc}{\ell^2}&1+\frac{2ac}{\ell^2}&0\\
0&0&0&0&1
\endpmatrix,
$$
where $\ell^2=2tb^2-2ac$.
$-\sigma_\ell\in \operatorname{SO }_\br^+$ if and only if $\ell^2>0$.
On the discriminant group $-\sigma_\ell$ defines
multiplication by
$$
\xi(\ell)=\frac{4t}{\ell^2}b^2-1.
$$
By definition $D=\hbox{div\,}(\ell)=(a,\,2tb,\,c)$, therefore
$D\,|\,2t$ and $D\,|\,\ell^2\,|\,2D$. We put
$$
\ell = {}^t(0,Da_1,\, b,\, Dc_1,0)
\qquad\text{with \ } (a_1, \,\frac{2tb}{D},\, c_1)=1.
$$
We have to consider four cases:
$$
\ell^2=D\quad\text{or}\quad \ell^2=2D\qquad\text{and}
\qquad D\, |\,t\quad\text{or}\quad D\hskip-2pt\not|  \ t.
$$

1). Let us suppose that  $D\,|\,t$ ($t_D=\dsize\frac{t}D$).
Then we have
$$\gather
\ell^2=2D\  \Longleftrightarrow \ 1=t_Db^2-Da_1c_1\
\Longleftrightarrow \ (D,t_D)=1\,\&\,\biggl( \frac{t_D}{D}\biggr)=1,\\
\xi(-\sigma_\ell)=2t_D b^2 -1 =
\cases -1&\operatorname{mod }\ 2t_D\\
       \hphantom{-} 1&\operatorname{mod }\ 2D
\endcases
\Longrightarrow -\sigma_\ell\in V_{t_D}\gt
\endgather
$$
(see \thetag{1.5}).
The case $\ell^2=D$ leads  trivially to a  contradiction.

2). Let us suppose that $D\,|\,2t$, but $\dsize\frac{t}D\not\in\bz$.
In this case
$D=2D_1$, $D_1|\,t$ and $t_{D_1}$ is odd. For such $D$ we have
$$
\gather
\ell^2=D\Longleftrightarrow 1=t_{D_1}b^2 - 4D_1 a_1c_1
\Longleftrightarrow \ (D_1,t_{D_1})=1\,\&\
(t_{D_1} \hbox{ odd})\,
\&\, \biggl( \dsize\frac{t_{D_1}}{4D_1}\biggr)=1,\\
\xi(-\sigma_\ell)=2t_{D_1} b^2 -1=
\cases -1&\operatorname{mod }\ 2t_{D_1}\\
      \hphantom{-}  1&\operatorname{mod }\ 2D_1
\endcases
\Longrightarrow -\sigma_\ell\in V_{t_{D_1}}\gt.
\endgather
$$
The case $\ell^2=2D$ leads to a contradiction to the  primitivity of
$\ell$.
\newline\qed\enddemo
For square free $t$ the system of the surfaces $\{\,H_{\ell^*}\}$,
where
$\ell^*=\dsize\frac{\ell}{\hbox{div\,}(\ell)}$
and  $\ell$ satisfies
the condition of   Lemma 3.8, contains all irreducible components
of the Humbert surfaces from Theorem 3.4.
This finishes the proof of Theorem 3.4.
\newline\qed

The next corollary follows immediately from  the proof of Lemma 3.8.

\proclaim{Corollary 3.9}Let $t$ be square free and $d||t$.
If $\biggl(\dsize\frac {d}{t_d}\biggr)=1$ and
$\biggl(\dsize\frac {d}{4t_d}\biggr)\ne 1$ then there is,
up to conjugation with respect to $\gt$,
exactly one, and if $\biggl(\dsize\frac {d}{4t_d}\biggr)=1$
then there are exactly two
involutions in $\Gamma_t V_d$.
They are
$\Psi^{-1}(-\sigma_{\ell_1})$ (in the  both cases),
and
$\Psi^{-1}(-\sigma_{\ell_2})$ (in the second case),
where
$$\align
\ell_1&={}^t(0,a_1,\dsize\frac {b_1}{t_d}, c_1,0)\in\widehat{L}_t,
\qquad
\,(a_1,b_1,c_1)=1,\quad db_1^2-t_da_1c_1=1\\
\ell_2&={}^t(0,a_2,\dsize\frac {b_2}{2t_d}, c_2,0)\in\widehat{L}_t,
\qquad
(a_2,b_2,c_2)=1,\quad db_2^2-4t_da_2c_2=1.
\endalign
$$
\endproclaim

\remark{Remarks}
1. It is possible to apply Lemma 3.8, which has been proved for any
integer $t$,
to classify the  divisorial part of the ramification locus
of the covering $\Cal A_t\to \Cal A_t^*$ for any integer $t$.

2. Using the same method one can construct divisors
 on a homogeneous domain of type IV of  any dimension $n$.
\endremark
\medskip

The ramification locus can also  have components of smaller
dimension.
The proof of Lemma 3.6 shows us that the orthogonal group
can contain a rotation in the positive definite subplane
of the lattice $L_t$.

\example{Example}(Brasch)
The following example is due to Brasch. It shows that in general the
ramification locus of the map ${\Cal A}_t \rightarrow {\Cal A}^*_t$
contains other components apart from the divisorial part described above.
Let $t\equiv 1 \operatorname{\,mod\,}4$. For an integer $f>0$ put
$$
c=-f^2 t-1,\quad  g=f^2.
$$
Then the matrix
$$
N=\pmatrix
-f\sqrt{t}   &   1/\sqrt{t}   &   0           &   f\sqrt{t}\\
c\sqrt{t}    &   0            &   f\sqrt{t}   &   f^2 t\sqrt{t}\\
c\sqrt{t}    &   0            &   f\sqrt{t}   &   -c\sqrt{t}\\
0            &   1/\sqrt{t}   &   -1/\sqrt{t} &   0
\endpmatrix
$$
is an element of $\Gamma_t^*$. One immediately checks that
$N^2=-E_4$. The fixed point set $\hbox{Fix } N$ is a curve
(cf. \cite{B, Hilfssatz 2.5.3}). Moreover the curve $\hbox{Fix } N$ is not
contained in the fixed point set $\hbox{Fix } I$ of an involution $I$ in
$\Gamma_t^*$.
This would namely imply that $IN=-NI$. Now using
the explicit form for $N$ which follows from \cite{B, Hilfssatz 2.8}, a
lengthy but straightforward calculation shows $c=-1$, a contradiction.
\endexample

As a next step we want to interprete the surfaces $H_t$, resp. $H_{4t}$ as
moduli spaces of abelian surfaces with real multiplication. It is well
known that there is a close connection between Hilbert modular surfaces
and Humbert surfaces \cite{F}, \cite{vdG, chapter IX}. Here we want to
determine precisely which Hilbert modular surfaces correspond to $H_t$,
resp. $H_{4t}$. Consider the ring $\frak{o}$ of integers in  the
number field $\bq(\sqrt{t})$ and recall that
$$
\frak{o} = \bz + \bz \omega,\quad \omega=\frac 1 2 (1+\sqrt{t})\quad
\hbox{if } t\equiv 1 \operatorname{mod} 4
$$
resp.
$$
\frak{o} = \bz + \bz \omega,\quad \omega=\sqrt{t} \quad
\hbox{if } t\not\equiv 1
 \operatorname{mod} 4.
$$
The {\it Hilbert modular group} $\operatorname{SL}_2(\frak{o})$ acts on
$\bh_1\times \bh_1$ by
$$
\pmatrix
\alpha   &   \beta\\
\gamma   &   \delta
\endpmatrix
(z_1, z_2) = \left(\frac{\alpha z_1 + \beta}{\gamma z_1 + \delta},
\ \frac{\alpha ' z_2 + \beta '}{\gamma ' z_2 + \delta '}\right)
$$
where $\ '$ denotes the Galois automorphism
$\ \sqrt{t} \mapsto - \sqrt{t}$.
\ The quotient space
\newline
$Y=\operatorname{SL}_2(\frak{o})\backslash \bh_1\times \bh_1$ is the
standard {\it Hilbert modular surface} associated to $\bq(\sqrt{t})$.
Let $\sigma$ be the involution which interchanges the two factors of
$\bh_1 \times \bh_1$, i.e.
$$
\sigma (z_1, z_2)=(z_2, z_1).
$$
Then the {\it symmetric} Hilbert modular group is
$$
\operatorname{SL}_2^{\sigma}(\frak{o})
=\operatorname{SL}_2 (\frak{o}) \cup
\sigma \operatorname{SL}_2  (\frak{o})
$$
and $Y^{\sigma} =\operatorname{SL}_2^{\sigma}
(\frak{o})\backslash \bh_1 \times
\bh_1$ is the corresponding symmetric Hilbert modular surface. We shall
first consider the Humbert surface $H_t$. In particular we assume that
$t\equiv 1 \operatorname{mod} 4$. To every point $(z_1, z_2) \in \bh_1
\times \bh_1$ one can associate the lattice
$$
\Lambda_{(z_1, z_2)}=\frak{o}\pmatrix
z_1\\ z_2
\endpmatrix+
\frak{o}\pmatrix
1\\1
\endpmatrix
= \bz
\pmatrix
z_1\\z_2
\endpmatrix+
\bz\pmatrix
\eta z_1\\ \eta ' z_2
\endpmatrix
+ \bz \pmatrix
-\eta ' / \sqrt{t}\\
\eta / \sqrt{t}
\endpmatrix
+ \bz
\pmatrix
 \sqrt{t}\\ -\sqrt{t}
\endpmatrix
$$
where $\eta = \sqrt{t} \omega$. The form
$$
E((x_1, x_2), (y_1, y_2))=\operatorname{Im} \left( \frac {x_1
\bar{y}_1}{\operatorname{Im} z_1} + \frac {x_2\bar{y}_2}{\operatorname{Im}
z_2}\right)
$$
defines a Riemann form for $\Lambda_t$. With respect to the basis given
above this is just the alternating form $W_t$. The torus
$$
A_{(z_1, z_2)}=\bc^2/ \Lambda_{(z_1, z_2)}
$$
is hence a $(1,t)$-polarized abelian surface with real multiplication in
$\frak{o}$. The Hilbert modular surface $Y$ is the moduli space of these
objects. The 2:1 cover $Y\rightarrow Y^{\sigma}$ identifies abelian
surfaces whose real multiplication  differs by the Galois conjugation. We
have a ``forgetful" map
$$
\Phi : Y^{\sigma} \rightarrow {\Cal A}_t.
$$
\proclaim{Theorem 3.10}
Assume $t\equiv 1 \operatorname{mod} 4$. The Humbert
surface $H_t$ is the image of the symmetric Hilbert modular surface
$Y^{\sigma}$ under the natural morphism $\Phi : Y^{\sigma}
\rightarrow {\Cal A}_t$ which is of degree 1 onto its image.
\endproclaim
Before giving the proof we turn to the Humbert surfaces $H_{4t}$.
Consider the ring
$$
\frak{o}_2 = \bz + \bz \sqrt{t}.
$$
Note that this is an order in $\frak{o}$ if $t\equiv 1
\operatorname{mod} 4$, whereas $
\frak{o}=\frak{o}_2$ if $t\not\equiv 1  \operatorname{mod} 4$. Let
$$
{\tilde \frak{o}}_2 = \frac 12 \bz + \frac 12 \sqrt{t} \bz,\quad {\tilde
\frak{o}}_2^{-1} = 2 \bz + 2\sqrt{t} \bz.
$$
The group
$$
\operatorname{SL}_2 (\frak{o}_2,{\tilde{\frak{o}}}_2)=
\{\pmatrix
a   &   b\\
c   &   d
\endpmatrix
|\  a, d \in \frak{o}_2,\  b \in {\tilde{\frak{o}}}_2,\  c \in
{\tilde\frak{o}}_2^{-1}, \ ad-bc=1 \}
$$
acts on $\bh_1 \times \bh_1$ as well as its symmetric counterpart
$$
\operatorname{SL}_2^{\sigma} (\frak{o}_2,
{\tilde{\frak{o}}}_2)=\operatorname{SL}_2 (\frak{o}_2,
{\tilde{\frak{o}}}_2) \cup \sigma \operatorname{SL}_2 (\frak{o}_2,
{\tilde{\frak{o}}}_2).
$$
Let
$$
{\widetilde Y}
=\operatorname{SL}_2(\frak{o}_2,
{\tilde{\frak{o}}}_2)\backslash \bh_1\times \bh_1,
\quad {\widetilde Y}^{\sigma} =
\operatorname{SL}^{\sigma}_2(\frak{o}_2,
{\tilde{\frak{o}}}_2)\backslash \bh_1\times \bh_1
$$
be the corresponding Hilbert modular surfaces.
Again the Riemann form $E$
induces a $(1,t)$-polarization on the tori
$$
A_{(z_1, z_2)}=\bc^2 / \Lambda_{(z_1, z_2)}
$$
where for $t\equiv 1\operatorname{mod} 4$
$$
\Lambda_{(z_1, z_2)}=\bz\pmatrix
z_1\\
z_2
\endpmatrix + \bz\pmatrix
2 \eta z_1\\
2 \eta ' z_2
\endpmatrix+ \bz \pmatrix
- \eta '/ \sqrt{t}\\
\eta/ \sqrt{t}
\endpmatrix + \bz  \pmatrix
\sqrt{t} /2\\
-\sqrt{t}/2
\endpmatrix
$$
resp.  $t\not\equiv 1 \operatorname{mod} 4$
$$
\Lambda_{(z_1, z_2)}=\bz\pmatrix
z_1\\
z_2
\endpmatrix+ \bz \pmatrix
\omega z_1\\
\omega ' z_2
\endpmatrix+ \bz \pmatrix
1/2\\
1/2
\endpmatrix+ \bz \pmatrix
\omega /2\\
\omega '/2
\endpmatrix\quad
$$
Hence ${\widetilde Y}$, resp. ${\widetilde Y}^{\sigma}$ are moduli spaces
of $(1,t)$-polarized abelian surfaces with real multiplication in
$\frak{o}_2$ and as before we have a canonical map
$$
{\widetilde{\Phi}} : {\widetilde Y}^{\sigma} \rightarrow {\Cal A}_t.
$$
\proclaim{Theorem 3.11}
The Humbert surface $H_{4t}$ is the image of the symmetric Hilbert modular
surface ${\widetilde Y}^{\sigma}$ under the natural map
$\widetilde{\Phi} : {\widetilde Y}^{\sigma}\rightarrow {\Cal A}_t$
which is of degree 1 onto its image.
\endproclaim
\demo{Proof of Theorems 3.10 and 3.11}
We shall treat the case of $H_{4t}$ and $t\equiv 1 \operatorname{mod }4$
in detail and then comment on the other cases. The proof is similar to
the proofs in \cite{HL \S 0}, cf. also \cite{F, Abschnitt 3}. Let
$$
R=\pmatrix
1   &   2\eta\\
1   &   2 \eta '
\endpmatrix
$$
and consider the map
$$
\align
{\widehat{\Phi}} :  \bh_1 \times \bh_1 &\rightarrow \bh_2\\
               (z_1, z_2) &\mapsto\  ^tR
\pmatrix
z_1   &   0\\
0   &   z_2
\endpmatrix  R.
\endalign
$$
Then
$$
\operatorname{Im }\widehat{\Phi} =
\{-(t^2-t)\tau_1+2t\tau_2-\tau_3=0\}={\Cal H}'_{4t}
$$
and modulo
$$
X=\pmatrix
1   &   0   &   0   &   0\\
t   &   1   &   0   &   0\\
0   &   0   &   1   &   -t\\
0   &   0   &   0   &   1
\endpmatrix
\in \Gamma_t
$$
this is equivalent to ${\Cal H}_{4t}=\{t \tau_1-\tau_3=0\}$. Let
$A_{\widehat{\Phi}(z_1, z_2)}$ be the abelian surface associated to the
period matrix
$$
\left( ^tR \pmatrix
z_1   &   0\\
0   &   z _2
\endpmatrix \
R \pmatrix
1   &   0\\
0   &   t
\endpmatrix  \right).
$$
Then $A_{(z_1, z_2)}$ and $A_{\widehat{\Phi}(z_1, z_2)}$ are
isomorphic as polarized abelian surfaces since
$$
^tR
\pmatrix
z_1   &   2\eta z_1   &   -\eta'/\sqrt{t}   &   \sqrt{t}/2\\
z_2   &   2\eta' z_2  &    \eta/\sqrt{t}    &   -\sqrt{t}/2
\endpmatrix=\left(
^tR \pmatrix
z_1   &   0\\
0     &   z _2
\endpmatrix\  R
\pmatrix
1   &   0\\
0   &   t
\endpmatrix
\right).
$$
Hence $\widehat{\Phi}$ is a lift of the map ${\widetilde{\Phi}}$.
Next we consider the homomorphism
$$
\Psi : \operatorname{SL}_2(\br) \times
\operatorname{SL}_2(\br)\longrightarrow \operatorname{Sp }_4(\br)
$$
$$
\Psi(
\pmatrix
a_1   &   b_1\\
c_1   &   d_1
\endpmatrix,\
\pmatrix
a_2   &   b_2\\
c_2   &   d_2
\endpmatrix)=
\pmatrix
^tR   &   0\\
0     &   R^{-1}
\endpmatrix
\pmatrix
d(a_1, a_2)   &   d(b_1, b_2)\\
d(c_1, c_2)    &   d(d_1, d_2)
\endpmatrix
\pmatrix
^tR   &   0\\
0     &   R
\endpmatrix
$$
where $d(a_1, a_2)=
\pmatrix
a_1   &   0\\
0     &   a_2
\endpmatrix
$, etc.
Via the embedding
$$
\align
\operatorname{SL}_2(\bq(\sqrt{t}))   &   \rightarrow
\operatorname{SL}_2(\br) \times \operatorname{SL}_2(\br)\\
\pmatrix
\alpha   &   \beta\\
\gamma   &   \delta
\endpmatrix  &   \mapsto
(\pmatrix
\alpha   &   \beta\\
\gamma   &   \delta
\endpmatrix,\
\pmatrix
\alpha'   &   \beta'\\
\gamma'   &   \delta'
\endpmatrix )
\endalign
$$
this also defines a homomorphism
$$
\widehat{\Psi} : \operatorname{SL}_2(\bq(\sqrt{t})) \rightarrow
\operatorname{Sp}_4(\br).
$$
For $J=
\left(\smallmatrix
0   &   1\\
1   &   0
\endsmallmatrix\right)$
we find that
$$
S=\pmatrix
^tRJ\ ^tR^{-1}   &   0\\
0               &   R^{-1} JR
\endpmatrix \in \Gamma_t.
$$
Clearly $S$ is an involution. Setting $\widehat{\Psi}(\sigma)=S$ we can
extend $\widehat{\Psi}$ to a homomorphism
$$
\widehat{\Psi}_{\sigma} : \operatorname{SL}_2^{\sigma}
(\bq(\sqrt{t}))\rightarrow \operatorname{Sp}_4(\br)
$$
and one checks easily that $\widehat{\Phi}$ is
$\widehat{\Psi}_{\sigma}$-equivariant. Let $G$, resp. $G_{\br}$ be the
stabilizer of ${\Cal H}'_{4t}$ in $\Gamma_t$, resp.
$\operatorname{Sp}_4(\br)$. As in \cite{HL, Lemma 0.8}, cf. also
\cite{F, Korollar 3.2.8} one shows that $G_{\br}$ is the group generated
by the image of $\Psi$ and $S$. The result follows if we can prove that
$G=\widehat{\Psi}_{\sigma} (\operatorname{SL}_2^{\sigma}(\frak{o}_2,
{\tilde\frak{o}}_2))$. For this it is now enough to prove the following.
\enddemo
\proclaim{Lemma 3.12}Let
$\pmatrix
a_i   &   b_i\\
c_i   &   d_i
\endpmatrix\in \operatorname{SL}_2(\br); \ i=1,2$ and assume that
$$
M=\pmatrix
^tR   &   0\\
0   &   R^{-1}
\endpmatrix
\pmatrix
d(a_1, a_2)   &   d(b_1, b_2)\\
d(c_1, c_2)   &   d(d_1, d_2)
\endpmatrix
\pmatrix
^tR^{-1}   &   0\\
0          &   R
\endpmatrix
\in \Gamma_t.
$$
Then $a_1, d_1 \in \frak{o}_2,\  b_1 \in {\tilde\frak{o}}_2,\  c_1 \in
{\tilde\frak{o}}_2^{-1}$ and $ a_2=a_1',\  b_2=b_1',\  c_2=c_1',\  d_2=d_2'$.
\endproclaim
\demo{Proof of the lemma}
We write
$M=
\pmatrix
A   &   B\\
C   &   D
\endpmatrix$.
Straightforward calculation gives
$$
A=\frac 1{2(\eta' -\eta)}
\pmatrix
2 \eta' a_1- 2\eta a_2   &   -a_1 + a_2\\
4 \eta \eta' (a_1 - a_2)   &   -2\eta a_1 + 2 \eta' a_2
\endpmatrix.
$$
{}From $A_{12} \in \bz$ we find $(a_1-a_2)/ (2\sqrt{t}) \in \bz$, i.e.
$$
a_2=a_1+2 n \sqrt{t},\quad n\in \bz.
$$
$A_{11}\in \bz$ gives $a_1 \omega' + a_2 \omega \in \bz$.
Hence
$$
a_1+2n\sqrt{t} \omega' = a_1(\omega + \omega') + 2n\sqrt{t} \omega'=
a_1\omega +a_2\omega' \in \bz.   \tag{3.2}
$$
I.e. $ a_1\in \bz + \bz\sqrt{t}=\frak{o}_2$. We can write $a_1=\alpha +
\beta \sqrt{t}$ with $ \alpha, \beta \in \bz$. By \thetag{3.2}
we have $a_1+n(\sqrt{t} + t) \in \bz$ and hence $\beta=-n$.
But then $a_2=a_1-2\beta
\sqrt{t} = \alpha-\beta\sqrt{t}=a_1'$. Note also that for $a_1\in
\frak{o}_2$ and
$a_2=a_1'$ one has $A\in \pmatrix
\bz   &   \bz\\
t\bz  &   \bz
\endpmatrix$.
Similarly we obtain
$$
B=\pmatrix
b_1 + b_2   &   2 \eta b_1 + 2\eta' b_2\\
2\eta b_1+ 2\eta' b_2   &   4 \eta^2 b_1 + 4 {\eta'}^{2} b_2
\endpmatrix.
$$
Since $B_{11} \in \bz$ we find that $b_2=n-b_1$ for some $n\in \bz$. Using
this and $B_{12}\in t \bz$
 one concludes that
$$
(b_1\omega - b_2\omega') \in \frac{\sqrt{t}}{2} \bz.
$$
Hence
$$
b_1-n\omega' = b_1(\omega+\omega')-n\omega'=b_1 \omega - b_2\omega' \in
\frac{\sqrt{t}}{2} \bz.                                  \tag{3.3}
$$
This shows $b_1, b_2 \in {\tilde \frak{o}}_2$. Writing $b_1=(\alpha +
\beta \sqrt{t})/2$ with $\alpha, \beta \in \bz$  one has from \thetag{3.3}
$$
\frac{\alpha}2 + \frac{\beta}2 \sqrt{t} - n\omega' \in \frac{\sqrt{t}}2
\bz
$$
i.e. $\alpha = n$. Hence $b_2=n-b_1=(\alpha-\beta \sqrt{t})/2=b_1'$.
Conversely if $ b_1\in {\tilde\frak{o}}_2$ and $b_2=b_1'$, then
$B\in
\pmatrix
\bz   &  t\bz\\
t\bz  &  t\bz
\endpmatrix$.
For $C$ one computes
$$
C=\frac{1}{4t}
\pmatrix
4{\eta'}^{2} c_1+4\eta^2 c_2   &   -2\eta' c_1 - 2\eta c_2\\
-2\eta' c_1-2\eta c_2        &   c_1 + c_2
\endpmatrix.
$$
Comparing this with the situation for B one finds that $c_1 \in {\tilde
\frak{o}}_2^{-1}, c_2=c_1'$.
Finally
$$
D=\frac{1}{2(\eta'-\eta)}
\pmatrix
2\eta' d_1 - 2\eta d_2   &   4\eta \eta' (d_1-d_2)\\
-d_1 + d_2               &   -2\eta d_1 + 2\eta' d_2
\endpmatrix.
$$
This case is analogous to A and one obtains $d_1 \in \frak{o}_2,\
d_2=d_1'$. This proves  Lemma 3.12.
\newline\qed\enddemo

The case of the Humbert surface $H_t$, $t\equiv 1 \operatorname{mod } 4$
can be treated in the same manner if we choose
$R=\pmatrix
1   &   \eta\\
1   &   \eta'
\endpmatrix$.
This is analogous to \cite{HL, \S 0} where the case
$t=5$ was done.
Finally we have to deal with $H_{4t}$ in case $t\not\equiv 1
\operatorname{mod } 4$. Here we can choose
$R=\pmatrix
1   &   \omega\\
1   &   \omega'
\endpmatrix, \omega = \sqrt{t}$.
In this case
$$
\operatorname{Im } \widehat{\Phi}=\{t\tau_1 - \tau_3=0\}={\Cal H}_{4t}
$$
and the above arguments go through essentially unchanged.
\newline\qed

\enddocument

\end